\documentclass[12pt]{amsart}
\title{Group schemes of square free order}
\author{V. Kumar Murty, Ying Zong}

\address{Department of Mathematics\\University of Toronto}
\email{murty@math.utoronto.ca, zongying@math.utoronto.ca}
\date{} % delete this line to display the current date

\begin{document}

\maketitle
%\tableofcontents

%\section{}
%\subsection{}

\section{Introduction}

\smallskip

Let $G$ be a finite flat finitely presented group scheme over a scheme $S$, with structural morphism $f: G\to S$. The $\mathcal{O}_S$-algebra $\mathcal{A}=f_*\mathcal{O}_G$, as $\mathcal{O}_S$-module, is locally free of finite type; its rank, a locally constant function on $S$, is the \emph{order}, or \emph{rank}, of $G$. We will say that $G$ is of \emph{prime order} (resp. \emph{square free order}) if its order at each point of $S$ is either $1$ or a prime number (resp. is square free, i.e., non-divisible by prime squares). 

\smallskip

This note is to exhibit how group schemes of prime order, whose story involving beautifully with such as Jacobi sums was told long ago by Oort and Tate \cite{oort_tate}, build up all group schemes of square free order, commutative or not, in a simple manner. 

\smallskip

\smallskip

{\bf Theorem 1.1.} --- \emph{Let $S$ be a scheme and let $G$ be a finite flat finitely presented $S$-group scheme of square free order. Then there exists a Hochschild extension of $S$-group schemes of square free order, 
\[1\to G'\to G\to G''\to 1,\] splitable by a finite \'{e}tale surjective base change, where $G''$ is finite \'{e}tale over $S$, and $G'$ a direct sum, as $G''$-modules, of commutative finite flat finitely presented $S$-group schemes of prime order. If
\[1\to G'_{\iota} \to G\to G''_{\iota}\to 1,\] $\iota=1, 2$, are two such extensions, there is a third such extension
\[1\to G'_3\to G\to G''_3\to 1\] with $G'_3=G'_1\cap G'_2$.}

\smallskip

\smallskip

In other words, at the expense of a finite \'{e}tale descent datum, every group scheme of square free order $G$ is a semi-direct product $G'G''$ and hence determined by how each summand of $G'$ of prime order is acted on by the finite \'{e}tale group scheme $G''$. We prove the theorem in 2.12 after recollecting a few facts. Recall that, by \emph{Hochschild extension} of group schemes, it is meant that 
\[1\to G'(T)\to G(T)\to G''(T)\to 1\] is exact for every $S$-scheme $T$. 

\smallskip

\smallskip

\section{Proof of Theorem 1.1.}

\smallskip

{\bf Lemma 2.1.} --- \emph{Let $S$ be a scheme of residue characteristics zero and let $G$ be a finite flat finitely presented $S$-group scheme. Then $G$ is \'{e}tale over $S$.}

\begin{proof} Since $G$ is flat and of finite presentation over $S$, it suffices to verify that $G\times_Ss$ is \'{e}tale over $k(s)$ for every $s\in S$ (EGA IV 17.6.2). That is, one can suppose that $S$ is the spectrum of a field of characteristic zero. By Cartier (SGA 3 VII B 3.3.1), every algebraic group over a field of characteristic zero is smooth, and so every \emph{finite} algebraic group over a field of characteristic zero is \'{e}tale.

\end{proof}

\smallskip

{\bf Lemma 2.2.} --- \emph{Let $S$ be a scheme and let $G$ be a finite \'{e}tale $S$-group scheme.}

\smallskip

\emph{Then, for each $d\in \Gamma(S, \mathbf{Z})$, the set $U$ of points $s$ of $S$ such that $G(\overline{s})$ has a unique subgroup of order $d(s)$, where $\overline{s}$ is the spectrum of an algebraic closure of $k(s)$, is open and closed in $S$. Over $U$, $G\times_SU$ has a unique finite \'{e}tale normal sub-$U$-group scheme of order $d$.}

\begin{proof} The question is local on $S$. One can assume $S$ affine and, by a ``passage \`{a} la limite'', noetherian and connected. If $U$ is not empty, let $\overline{s}\to S$ be a geometric point with image $s$ in $U$ and let \[\rho: \pi_1(S, \overline{s})\to \mathrm{Aut}(G(\overline{s}))\] be the monodromy representation corresponding to the finite \'{e}tale $S$-group scheme $G$. The unique subgroup of $G(\overline{s})$ of order $d(s)$ is characteristic in $G(\overline{s})$ and in particular invariant by $\rho$. So $G$ has a unique finite \'{e}tale normal sub-$S$-group scheme of order $d$, and $U=S$.

\end{proof}

\smallskip

{\bf Lemma 2.3.} --- \emph{Let $S$ be the spectrum of an algebraically closed field of characteristic $p>0$ and let $G$ be a finite $S$-group scheme of square free order. Then, the identity component of $G$, if non-trivial, is $\mu_p$ or $\alpha_p$.}

\begin{proof} This follows by SGA 3 XVII 4.2.1.

\end{proof}

\smallskip

{\bf Lemma 2.4.} --- \emph{Let $S$ be a scheme and let $G$ be a finite flat finitely presented $S$-group scheme of square free order. If, for every point $s$ of $S$, $i_G(s)$ denotes the order of the identity component of $G\times_Ss$, the function $s\mapsto i_G(s)$, called the infinitesimal rank of $G$, is locally constructible and upper semi-continuous.}

\begin{proof} In fact, the order of the quotient of $G\times_Ss$ by its identity component, the \emph{separable rank} of $G\times_Ss$, is a locally constructible lower semi-continuous function of $s\in S$ (EGA IV 15.5.1).

\end{proof}

\smallskip

{\bf Lemma 2.5.} --- \emph{Keep the notations and hypotheses of $2.4$. The set $S_1$ of points $s$ of $S$ such that $i_G(s)$ is $1$ is a retro-compact open subset of $S$; it consists exactly of all those points of $S$ above which $G$ is $S$-\'{e}tale.}

\begin{proof} This follows by 2.4 and EGA IV 17.6.2.

\end{proof}

\smallskip

{\bf Lemma 2.6.} --- \emph{Let $S$ be a scheme and let $G$ be a finite flat finitely presented $S$-group scheme of prime order. Then the group law of $G$ is commutative.}

\begin{proof} This is \cite{oort_tate} \S 1 Theorem 1.

\end{proof}

\smallskip

{\bf Lemma 2.7.} --- \emph{Let $S$ be a scheme and let $G$ be a commutative finite flat finitely presented $S$-group scheme. Then $G$ is of $n$-torsion, $n\in \Gamma(S, \mathbf{N})$ its order, and hence a module over the locally constant sheaf of rings $\mathbf{Z}/n\mathbf{Z}$.}

\begin{proof} This is \cite{oort_tate} \S 1 Theorem.

\end{proof}

\smallskip

{\bf Lemma 2.8.} --- \emph{Let $S$ be a scheme, $n=\prod p$ a prime factorization of a square free section $n\in\Gamma(S, \mathbf{N})$, and $G$ a commutative finite flat finitely presented $S$-group scheme of order $n$.}

\smallskip

\emph{Then $G$ is in a unique way a direct sum, indexed by these prime factors $p\in \Gamma(S, \mathbf{N})$ of $n$, of finite flat finitely presented $S$-group schemes of order $p$. This decomposition commutes with every base change and is invariant by all $S$-group automorphisms of $G$.}

\begin{proof} As a module over 
\[\mathbf{Z}/n\mathbf{Z}=\prod \mathbf{Z}/p\mathbf{Z}\] (2.7), $G$ is the direct sum of its $p$-torsion sub-modules $G_p$; this decomposition is clearly invariant by all $S$-group automorphisms of $G$ and compatible with base change.
\smallskip

Being the kernel of $p.\mathrm{Id}_G$, $G_p$ is a finitely presented closed sub-$S$-group scheme of $G$. And $G_p$ is also the \emph{fppf} image of $(n/p).\mathrm{Id}_G$, as one verifies fiberwise (EGA IV 11.3.11) and uses 2.1+2.3; consequently, $G_p$ is $S$-flat (EGA IV 11.3.11). The order $n_p$ of $G_p$, which is divisible by $p$ as is seen fiber by fiber, is $p$, since
\[n=\prod p\ |\ \prod n_p=n.\] 

Any decomposition of $G$ into group schemes of prime order $p$ necessarily coincides with the the above, as one sees immediately by 2.6+2.7.

\end{proof}

\smallskip

{\bf Lemma 2.9.} --- \emph{Let $S$ be a scheme and let 
\[1\to G'\to G\to G''\to 1\] be an extension of a finite \'{e}tale $S$-group scheme $G''$ by a commutative finite flat finitely presented $S$-group scheme $G'$. Assume that $G$ is of square free order.}

\smallskip

\emph{Then this extension is Hochschild and splitable by a finite \'{e}tale surjective base change.}

\begin{proof} One may evidently assume $G''$ of constant order $n''$.  
\smallskip

Notice that $n''.\mathrm{Id}_{G'}: G'\to G'$ is an isomorphism. Namely, $n''.\mathrm{Id}_{G'\times_Ss}$ is an automorphism of $G'\times_Ss$ for every $s\in S$. Indeed, as $G$ is of square free order, $G'$ is of order coprime to $n''$ at every point $s\in S$.

\smallskip

Let $T$ be an arbitrary $S$-scheme. Consider the cohomology sequence
\[1\to G'(T)\to G(T)\to G''(T)\stackrel{d}{\to} H^1((\mathrm{Sch}/T)_{\mathrm{fppf}}, G')\] For each $x\in G''(T)$, one has $d(x)=0$, since
\[n''.d(x)=d(x^{n''})=d(1)=0\] and $n''.\mathrm{Id}_{G'}: G'\to G'$ is an isomorphism. Thus the extension
\[1\to G'\to G\to G''\to 1\] is Hochschild. 
\smallskip

For the second statement, at the expense of a finite \'{e}tale surjective base change, one may assume that the finite \'{e}tale $S$-group scheme $G''$ is locally constant and then constant of value $|G''|$. It suffices to show that the exact sequence
\[1\to G'(S)\to G(S)\to G''(S)\to 1\] is split by $|G''|\to G''(S)$. This follows as $H^2(B|G''|, G'(S))=0$.

\end{proof}

\smallskip

{\bf Lemma 2.10.} --- \emph{Let $S$ be a scheme and let $G$ be a finite flat finitely presented $S$-group scheme of square free order. If, for $\iota=1, 2$, either
\[1\to G'_{\iota}\to G\to G''_{\iota}\to 1\] writes $G$ as extension of a finite \'{e}tale $S$-group scheme by a finite flat finitely presented $S$-group scheme, then there is a third such extension 
\[1\to G'\to G\to G''\to 1\] with $G'=G'_1\cap G'_2$.}

\begin{proof} The intersection $G'=G'_1\cap G'_2$ is a finitely presented normal closed sub-$S$-group scheme of $G$. It remains to show that $G'$ is flat over $S$ and that the quotient $G''=G/G'$, which is then representable by a finite flat finitely presented $S$-group scheme, is \'{e}tale over $S$. 
\smallskip

The assertion being local on $S$, one may assume $S$ affine, then by a ``passage \`{a} la limite'', noetherian and local, then by completion along its closed point, complete, and finally by reduction to each of its closed sub-schemes of finite lengths, artin local. Let the closed point of $S$ be $s$ and let $G_s^o$ be the identity component of $G_s=G\times_Ss$. Consider the canonical exact sequence
\[1\to G_s^o\to G_s\to G_s/G_s^o\to 1.\] 
As the base change by $s\to S$ induces an equivalence from the category of finite \'{e}tale $S$-group schemes onto the category of finite \'{e}tale $s$-group schemes (EGA IV 18.1.2), there is up to unique isomorphisms a unique finite \'{e}tale $S$-group scheme $Q$ with $Q\times_Ss=G_s/G_s^o$. By EGA IV 18.1.3, there is a unique $S$-group homomorphism $G\to Q$ which specializes to the projection $G_s\to G_s/G_s^o$. This morphism $G\to Q$ is faithfully flat (EGA IV 11.3.11) and provides an exact sequence of $S$-group schemes 
\[1\to P\to G\to Q\to 1\] with $P$ finite flat over $S$, $P\times_Ss=G_s^o$. For $\iota=1, 2$, either composition
\[P\to G\to G''_{\iota}\] is by EGA IV 18.1.3 trivial, for it is when specialized to $s$. So $P\subset G'_1\cap G'_2=G'$. Replacing $G$ by $G/P$, $G'_{\iota}$ by $G'_{\iota}/P$ and $G'$ by $G'/P$, one can assume $G=Q$ finite \'{e}tale over $S$. But now all of $G'_1, G'_2, G', G''$ are finite \'{e}tale over $S$. This completes the proof.

\end{proof}

\smallskip

{\bf Lemma 2.11.} --- \emph{Keep the notations and hypotheses of $2.4$ and let $p$ be a prime number.}

\smallskip

\emph{Then the set $S_p$ of points $s$ of $S$ such that $i_G(s)$ divides $p$ is a retro-compact open subset of $S$. The union $V_p$ of $S_p-S_1$ and the points $s$ of $S$ such that $G(\overline{s})$ has a unique $p$-Sylow subgroup of order $p$, where $\overline{s}$ is the spectrum of an algebraic closure of $k(s)$, is open and closed in $S$. On $V_p$, $G\times_SV_p$ has a unique finite flat finitely presented normal sub-$V_p$-group scheme $G_p$ of order $p$.}

\begin{proof} That $S_p$ is a retro-compact open subset of $S$ follows from 2.4+2.5. 

\smallskip

I) \emph{The set $U_p:=V_p\cap S_p$ is open and closed in $S_p$, and $G\times_SU_p$ has a unique finite flat finitely presented normal sub-$U_p$-group scheme $P$ of order $p$. The quotient $(G\times_SU_p)/P$ is finite \'{e}tale over $U_p$ }:

\smallskip

Restricting to $S_p$, one may suppose $S=S_p$, namely, that $i_G$ takes values in $\{1, p\}$. Let $U: (\mathrm{Sch}/S)^o\to (\mathrm{Sets})$ be the following sub-functor of the final functor :
\smallskip

\emph{For an $S$-scheme $S'$, $U(S')=\{\emptyset\}$, if $G\times_SS'$ is the extension of a finite \'{e}tale $S'$-group scheme by a finite flat finitely presented $S'$-group scheme of order $p$, and $U(S')=\emptyset$, otherwise.}

\smallskip

It suffices to show that $U\to S$ is representable by an open and closed immersion. 
\smallskip

\emph{The functor $U$ is locally of finite presentation }:

\smallskip

This is immediate by the technique of ``passage \`{a} la limite''.

\smallskip

\emph{The functor $U$ is a sheaf on $\mathrm{Sch}/S$ for the \'{e}tale topology }:

\smallskip

Namely, $U(T)=U(T')$ for every \'{e}tale surjective morphism $T'\to T$ of $S$-schemes. Indeed, if $U(T')$ is non-empty so that there is an extension of $T'$-group schemes
\[1\to P'\to G\times_ST'\to Q'\to 1\] with $Q'$ $T'$-finite \'{e}tale and $P'$ finite flat of finite presentation over $T'$ of order $p$, it follows by 2.10 that $p_1^*P'=p_2^*P'$ in $G\times_ST''$, where $p_1, p_2$ are the two projections of $T''=T'\times_TT'$ onto $T'$. On $P'\hookrightarrow G\times_ST'$, there exists therefore a canonical descent datum relative to $T'\to T$, and so $U(T)$ is non-empty as well. 

\smallskip

\emph{The functor $U$ is formally \'{e}tale }:
\smallskip

Namely, $U(T)=U(T')$ for every nilpotent $S$-immersion $T'\hookrightarrow T$. This is immediate from EGA IV 18.1.2+18.1.3 by a similar argument as in 2.10.

\smallskip

\emph{The functor $U$ verifies the valuative criterion of properness }:
\smallskip

As $U$ is a sheaf on $\mathrm{Sch}/S$ for the \'{e}tale topology, the question is local on $S$. It suffices to assume $S$ affine and, by a ``passage \`{a} la limite'', noetherian. Given an $S$-scheme $T$ which is the spectrum of a discrete valuation ring and which has generic point $t$, we need to prove that $U(T)=U(t)$. In fact, if $U(t)$ is non-empty and so $G\times_St$ is the extension of a finite \'{e}tale $t$-group scheme $Q_t$ by a finite $t$-group scheme $P_t$ of order $p$, then $Q=G_T/P$ is finite \'{e}tale over $T$, where $P$ denotes the closed image of $P_t$ in $G_T=G\times_ST$. For, the infinitesimal rank of $G$, hence that of $Q$ as well, takes values in $\{1, p\}$, and so $i_{Q}(T)=\{1\}$, as $Q$ is of order prime to $p$.

\smallskip

It is now clear that $U$ is representable by an open and closed sub-scheme of $S$ with underlying set $U_p$.

\smallskip

\smallskip

II) \emph{The set $U'_p:=V_p-(S_p-S_1)$ is open and closed in $S-(S_p-S_1)$, and $G\times_SU'_p$ has a unique finite \'{e}tale normal sub-$U'_p$-group scheme of order $p$ }:

\smallskip

Restricting to each $S_q$ at a time, where $q$ is a prime distinct from $p$, one may assume $S=S_q$. By 2.2 and by I), one may assume $S=U_q$. Consider the canonical exact sequence provided in I)
\[1\to Q\to G\to G/Q\to 1\] where $Q$ is the unique finite flat finitely presented normal sub-$S$-group scheme of $G$ of order $q$. By 2.2 applied to the prime-to-$q$ finite \'{e}tale quotient $G/Q$, and replacing $S$ by an open and closed sub-scheme, one may assume $G/Q$ has a unique finite \'{e}tale normal sub-$S$-group scheme $\overline{P}$ of order $p$. Let $\overline{G}$ be the pre-image of $\overline{P}$ in $G$. For proving II), one may replace $G$ by $\overline{G}$ and assume $P:=G/Q$ of order $p$. By the uniqueness assertion of II) and according to 2.9, at the expense of a finite \'{e}tale descent datum, one may assume that $G=QP$ is a semi-direct product and that $P$ is constant of generator $g$. Now, by conjugation, $P$ acts on $Q$. Let $\zeta$ be the image of $g$ in $\mathrm{Aut}_S(Q)$ and we will see that the following sub-functor of the final functor is open and closed :
\smallskip

\emph{For an $S$-scheme $S'$, $V(S')=\{\emptyset\}$, if $\zeta\times_SS'=\mathrm{id}_Q\times_SS'$, and $V(S')=\emptyset$, otherwise.}

\smallskip

This sub-functor $V\to S$ is representable by a closed immersion (SGA 3 VIII 6.4). Clearly, $V$ has underlying set $U'_p$, and $G\times_SV$ is commutative with its $p$-torsion being the unique normal finite \'{e}tale sub-$V$-group scheme of order $p$. 

\smallskip

On $S_1$, $V\times_SS_1\to S_1$ is an open and closed immersion. For, $Q\times_SS_1$, hence $\underline{\mathrm{Aut}}_S(Q)\times_SS_1$ as well, is finite \'{e}tale over $S_1$.

\smallskip

To show that $V\to S$ is an open immersion at each point $s$ of $U'_p-S_1$, it suffices by strictly localizing $S$ at $s$ to assume $S$ strictly local with closed point of characteristic $q$. By \cite{oort_tate} \S 2 Theorem 2, $Q$ corresponds to a triple $(L, c, d)$ which consists of an invertible $\mathcal{O}_S$-module $L$ and of $\mathcal{O}_S$-linear homomorphisms $c: L\to L^{\otimes q}$, $d: L^{\otimes q}\to L$ satisfying $d\circ c=w.\mathrm{Id}_L$ for a certain element $w\in\Gamma(S, \mathcal{O}_S)$. This correspondence identifies the $S$-group automorphisms of $Q$ with the units $u\in \Gamma(S, \mathcal{O}_S)^{\times}$ such that $uc=u^qc$, $ud=u^qd$. In particular, $\zeta$ is identified to an element of $\mu_p(S)$, and the relation ``$\zeta=\mathrm{Id}_Q$'' is an open and closed relation on $S$, since $\mu_{p S}$ is finite \'{e}tale over $S$.

\smallskip

The combination of I) and II) finishes the proof.

\end{proof}

\smallskip

2.12. \emph{Proof of Theorem 1.1.}

\smallskip

The last assertion follows by 2.10. By 2.8+2.9, it only remains to write $G$ as an extension
\[1\to G'\to G\to G''\to 1\] with $G''$ $S$-finite \'{e}tale and $G'$ commutative finite flat of finite presentation over $S$. 
\smallskip

One may assume $G$ of constant order $n$. The order $i_G(s)$ of the identity component of each fiber $G\times_Ss$ divides $n$ and is equal to either $1$ or the characteristic of $k(s)$ (2.3). If $E:=i_G(S)$ consists only of $1$, $G$ is $S$-finite \'{e}tale (2.5). Otherwise, for each prime $p$ in $E$, there exists (2.11) an open and closed sub-scheme $V_p$ of $S$ such that :
\smallskip

--- \emph{On $S-V_p$, $i_G$ does not take value $p$.}

\smallskip

--- \emph{On $V_p$, $G\times_SV_p$ has a unique finite flat finitely presented normal sub-$V_p$-group scheme $G_p$ of order $p$.}

\smallskip

\smallskip

By induction on the cardinality of $E-\{1\}=\{p_1, p_2, \cdots\}$, suppose that the claim holds already for the restriction of $G$ to 
\[S-(V_{p_1}\cap V_{p_2}\cap\cdots)=(S-V_{p_1})\coprod (V_{p_1}-V_{p_2})\coprod (V_{p_1}\cap V_{p_2}-V_{p_3})\coprod \cdots\] Then, restricting to $V=V_{p_1}\cap V_{p_2}\cap\cdots$, one may assume $V_p=S$ for every prime $p$ in $E$. Put $G'=G_{p_1}\times_SG_{p_2}\times_S\cdots$. Now, the $S$-morphism
\[G'\to G,\ (g_1, g_2, \cdots)\mapsto g_1.g_2\cdots\] identifies $G'$ as a normal sub-$S$-group scheme of $G$ with $G''=G/G'$ finite \'{e}tale over $S$.

\smallskip

\smallskip

\smallskip

\smallskip

%-----------------------------------------------------------------------------
\bibliographystyle{amsplain}

\end{document}